\documentclass[11pt, a4paper]{article}
\usepackage[latin1]{inputenc}
\usepackage{amsfonts}

\usepackage{theorem}

\newcommand{\R}{\mathbb{R}}

\newcommand{\qed}{\hfill  $\Box$ }

\theoremstyle{plain}
\theorembodyfont{\rmfamily}
\newtheorem{defn}{Definition}[section]

\theorembodyfont{}
\newtheorem{lem}[defn]{Lemma}

\newtheorem{thm}[defn]{Theorem}
\newtheorem{cor}[defn]{Corollary}

\begin{document}

\title{Curvature estimates for graphs with prescribed mean curvature and flat normal bundle}

\author{Steffen Fröhlich, Sven Winklmann} 

\date{\today}

\maketitle

\begin{abstract}

We consider graphs $\Sigma^n \subset \R^m$ with prescribed mean curvature and flat normal bundle. Using techniques of Schoen, Simon and Yau \cite{SSY} and Ecker-Huisken \cite{Ecker:Huisken:HPC}, we derive the interior curvature estimate 
$$\sup_{\Sigma \cap B_R} |A|^2 \leq \frac{C}{R^2}$$
up to dimension $n \leq 5$, where $C$ is a constant depending on natural geometric data of $\Sigma$ only. This generalizes previous results of Smoczyk, Wang and Xin \cite{SWX} and Wang \cite{Wang:stability} for minimal graphs with flat normal bundle.\\ 


\end{abstract}

\section{Introduction}

Let $\psi :\Omega \rightarrow \R^k$ be a smooth function defined on a domain $\Omega \subset \R^n$, and denote by $\Sigma =\{(x,\psi(x)) : x \in \Omega \}$ the corresponding graph in $\R^{m=n+k}$. In this paper we assume the normal bundle of $\Sigma$ to be flat and prove the interior curvature estimate   
\begin{equation} \label{main:sketchy}
\sup_{\Sigma \cap B_R} |A|^2 \leq \frac{C}{R^2}
\end{equation}
up to dimension $n \leq 5$, where $|A|$ denotes the length of the second fundamental form, $B_R \subset \R^m$ is a closed ball of radius $R$ centered at some point $p \in \Sigma$, and $C$ is a constant depending on natural geometric data of $\Sigma$ only, see Theorem \ref{main:thm}. 

Recently, curvature estimates for minimal graphs with flat normal bundle have been established independently by Smoczyk, Wang and Xin \cite{SWX} and Wang \cite{Wang:stability}. In particular, they have obtained higher dimensional analogues of the famous Schoen-Simon-Yau estimates \cite{SSY} and Ecker-Huisken's Bernstein result \cite{Ecker:Huisken} for entire minimal graphs of controlled growth. 

Without any geometric restrictions on the normal bundle the situation turns out to be more complicated as can be seen from the counter example of Lawson-Osserman \cite{Lawson:Osserman}. In \cite{HJW} Hildebrandt, Jost and Widman have studied entire solutions of the minimal surface system
$$\frac{\partial}{\partial x^i} \left( \sqrt{g} g^{ij} \frac{\partial \psi^\alpha}{\partial x^j} \right) = 0, \quad \alpha=1, \ldots, k.$$
Here, $g_{ij}=\delta_{ij} + \sum_\alpha \frac{\partial \psi^\alpha}{\partial x^i} \frac{\partial \psi^\alpha}{\partial x^j}$, $(g^{ij})=(g_{ij})^{-1}$ and $g=\det(g_{ij})$. Using a regularity estimate for harmonic maps they could prove a Bernstein result under a suitable lower bound on the function 
$$w = \left[ \det \left( \delta_{ij} + \sum_\alpha D_i \psi^\alpha D_j \psi^\alpha \right) \right]^{-1/2}.$$  
Later, their result has been improved by Jost-Xin \cite{Jost:Xin} and Wang \cite{Wang:Bernstein}. In fact, Wang's Bernstein result holds for the entire class of area decreasing maps with bounded gradient. For a detailed survey on minimal graphs in higher co-dimension and further comments on the literature we refer to the recent monograph of Giaquinta-Martinazzi \cite[Chapter 11]{GM}. We also remark, that more explicit estimates for two-surfaces in $\R^m$ can been obtained by using strictly two-dimensional techniques, cf. Osserman \cite{Osserman} and Bergner-Fröhlich \cite{Bergner:Froehlich}. 

The paper is organized as follows: In section 2 we first collect some basic facts on graphs with flat normal bundle. Using ideas of Ecker-Huisken \cite{Ecker:Huisken:HPC} we then prove a rather general Simons inequality (Lemma \ref{Simons:inequality:flat}) for the Laplacian of the length of the second fundamental form. In section 3 we use this estimate to derive the $L^p$ curvature bound 
$$\int_{\Sigma \cap B_R} |A|^p \, d\mathcal{H}^n \leq C R^{n-p}$$
for some $p>n$ with a constant $C$ depending only on the geometric data of the problem, see Theorem \ref{thm:integral:curvature:estimate}. Here, we can proceed similarly as Winkl\-mann \cite{winklmann:manuscripta} who established a corresponding estimate for hypersurfaces of prescribed anisotropic mean curvature. In view of a general mean value inequality (Lemma \ref{meanvalue}), which is of independent interest on its own, this leads to the desired curvature estimate (\ref{main:sketchy}). As an application of our results we recover the Bernstein result of Smoczyk, Wang and Xin \cite{SWX} and Wang \cite{Wang:stability} for minimal graphs with flat normal bundle. \\

\noindent \textbf{Acknowledgement.} The second author was financially supported by the Alexander von Humboldt foundation and the Centro di Ricerca Matematica Ennio De Giorgi via a Feodor Lynen research scholarship.

\section{Notation and preliminary results}

Let $f: \Sigma^n \rightarrow \R^{m=n+k}$ be a smooth immersion of an $n$-dimensional, oriented manifold without boundary into euclidean $m$-space of arbitrary co-dimension $k \geq 1$. We denote by
$$g(X,Y) = \langle df(X), df(Y) \rangle $$
the induced metric with corresponding Levi-Civita connection
$$\nabla_X Y = (D_X Y)^\top $$
and curvature tensor
$$R(X,Y)Z = \nabla_X \nabla_Y Z - \nabla_Y \nabla_X Z - \nabla_{[X,Y]} Z.$$
Here, $X,Y,Z$ are smooth vectorfields on $\Sigma$, $D$ denotes the covariant derivative on $\R^m$ and $(\cdot)^\top$ is the projection onto $T \Sigma$, the tangent bundle of $\Sigma$, which we will always identify with $df(T \Sigma)$. 

The second fundamental form is given by
$$A(X,Y) = (D_X Y)^\perp = D_X Y - \nabla_X Y,$$
where $(\cdot)^\perp$ is the projection onto the normal bundle $N \Sigma$. Taking its trace defines the mean curvature vector
$$H=\mbox{trace}(A).$$

We also have an induced connection on the normal bundle $N \Sigma$ defined by the relation
$$\nabla^\perp_X \eta = (D_X \eta)^\perp$$
for any normal section $\eta$. The corresponding curvature tensor is given by
$$R^\perp(X,Y)\zeta = \nabla^\perp_X \nabla^\perp_Y \zeta - \nabla^\perp_Y \nabla^\perp_X \zeta - \nabla^\perp_{[X,Y]} \zeta.$$

We remark that these connections extend naturally to higher order tensor bundles formed from $T \Sigma$ and $N \Sigma$. For example, for an $(0,r)$-tensor $T$ with values in $N \Sigma$ the covariant derivative $\nabla T$ is given by
\begin{eqnarray*}
(\nabla_X T)(Y_1 , \ldots, Y_r) & = & \nabla^\perp_X T(Y_1, \ldots, Y_r) - T(\nabla_X Y_1, \ldots, Y_r) \\
& & - \ldots - T(Y_1, \ldots , \nabla_X Y_r).
\end{eqnarray*}

Let us now choose local orthonormal frames $\{e_i\}_{i=1, \ldots, n}$ and $\{e_\alpha\}_{\alpha=n+1, \ldots, m}$ for $T \Sigma$ and $N \Sigma$, respectively. In these frames the coefficients of the second fundamental form are given by 
$$h_{\alpha i j} = \langle A(e_i, e_j) , e_\alpha \rangle = - \langle D_{e_i} e_\alpha, e_j \rangle$$
and the mean curvature vector by
$$H=H_\alpha e_\alpha \quad \mbox{with } H_\alpha = h_{\alpha i i}.$$
Here and in the following we are using Einstein's summation convention: Repeated Latin and Greek indices are automatically summed from $1$ to $n$ and from $n+1$ to $m$, respectively, unless not otherwise stated. We also write
$$\nabla_k h_{\alpha i j} = \langle (\nabla_{e_k} A)(e_i, e_j) , e_\alpha \rangle,$$
$$R_{ijkl} = g( R(e_i, e_j) e_k, e_l )$$ 
and 
$$R^\perp_{ij\alpha \beta} = \langle R(e_i, e_j) e_\alpha, e_\beta \rangle$$
for the coefficients of $\nabla A$, $R$ and $R^\perp$. The fundamental equations of Gauß, Codazzi and Ricci then take the form 
\begin{equation} \label{Gauss}
R_{ijkl} = h_{\alpha il} h_{\alpha jk} - h_{\alpha ik} h_{\alpha jl},
\end{equation}
\begin{equation} \label{Codazzi}
\nabla_k h_{\alpha ij} = \nabla_i h_{\alpha jk}
\end{equation}
and
\begin{equation} \label{Ricci}
R^\perp _{ij \alpha \beta} = h_{\beta ik} h_{\alpha jk} - h_{\alpha ik} h_{\beta jk}.
\end{equation}

We also write $\nabla_i \nabla_j \varphi$ for the coefficients of $\nabla \nabla \varphi$, the second covariant derivative of a smooth function $\varphi$. The Laplace-Beltrami operator is then given by $\Delta \varphi = \nabla_i \nabla_i \varphi$. More generally, for any $(0,r)$-tensor with values in $N\Sigma$ we write $\nabla_i \nabla_j T_{\alpha k_1 \ldots k_r}$ for the coefficients of $\nabla \nabla T$. Finally, we denote by $|T|^2 = \sum_{i_1, \ldots, i_r} |T(e_{i_1}, \ldots, e_{i_r})|^2$ the square of the length of $T$.  

The following identity was first proved by Simons \cite{Simons} and is a direct consequence of (\ref{Gauss}), (\ref{Codazzi}) and (\ref{Ricci}). For further details see also Wang \cite[Section 7]{Wang:mcf}. 

\begin{lem}
For an arbitrary immersion $f: \Sigma^n \rightarrow \R^m$ the second fundamental form satisfies 
\begin{eqnarray} \label{Simons}
\frac{1}{2} \Delta |A|^2 & = & |\nabla A|^2 +  h_{\alpha ij} \nabla_i \nabla_j H_\alpha  + H_\alpha h_{\alpha ij} h_{\beta jk} h_{\beta ki} \nonumber \\
& &  - \sum_{i,j,k,l} (  h_{\alpha ij} h_{\alpha kl} )^2 - |R^\perp|^2. 
\end{eqnarray}
\end{lem}

Next, we consider the parallel $n$-form $\Omega = dx^1 \wedge \ldots \wedge dx^n$ on $\R^m$ and put
$$w = *\Omega = \Omega(e_1, \ldots, e_n),$$
where $*$ is the Hodge operator. The following equation is due to Fischer-Colbrie \cite{FC} and Wang \cite{Wang:staromega}, \cite{Wang:Bernstein}. For an alternative exposition we also refer to Giaquinta-Martinazzi \cite[Chapter 11]{GM}:

\begin{lem}
For an arbitrary immersion $f: \Sigma^n \rightarrow \R^m$ the function $w=*\Omega$ satisfies
\begin{equation} \label{Jacobi}
\Delta w + |A|^2 w = \Omega_{\alpha i} \nabla_i H_\alpha - 2 \sum_{\alpha < \beta, i<j} \Omega_{\alpha \beta ij} R^\perp_{ij\alpha\beta},  
\end{equation}
where $\Omega_{\alpha i} = \Omega(e_1, \ldots, e_\alpha, \ldots, e_n)$ with $e_\alpha$ occupying the $i$-th position, and $\Omega_{\alpha \beta ij} = \Omega(e_1, \ldots, e_\alpha, \ldots, e_\beta, \ldots , e_n)$ with $e_\alpha$, $e_\beta$ occupying the $i$-th and $j$-th position, respectively.
\end{lem}

In this paper we are particularly interested in immersions with flat normal bundle, that is the case $R^\perp =0$. The above equations then simplify as follows:
\begin{eqnarray} \label{Simons_flat}
\frac{1}{2} \Delta |A|^2  & = & |\nabla A|^2 +  h_{\alpha ij} \nabla_i \nabla_j H_\alpha  + H_\alpha h_{\alpha ij} h_{\beta jk} h_{\beta ki}  \nonumber \\
& & - \sum_{i,j,k,l} (  h_{\alpha ij} h_{\alpha kl} )^2
\end{eqnarray}
and
\begin{equation} \label{Jacobi_flat}
\Delta w + |A|^2 w =  \Omega_{\alpha i} \nabla_i H_\alpha .  
\end{equation}

Suppose now that $\Sigma = \{(x,\psi(x)) : x \in \Omega \}$ is the graph of a smooth function $\psi: \Omega \rightarrow \R^k$ over some domain $\Omega \subset \R^n$. In this case one easily checks the identity
$$w= \left[ \det \left( \delta_{ij} + D_i \psi^\alpha D_j \psi^\alpha \right) \right]^{-1/2}.$$
In particular we have $w>0$. Define the quantity $K_1$ by
\begin{equation} \label{Definition_H1}
K_1 = \left(w^{-1} \Omega_{\alpha i} \nabla_i H_\alpha \right)^+, 
\end{equation}
where $g^+$ denotes the positive part of the function $g$. Moreover, denote by $\mathcal{H}^n$ the $n$-dimensional Hausdorff measure. Then we can state an energy-type estimate as follows:

\begin{lem} 
Suppose $\Sigma^n \subset \R^m$ is a graph with flat normal bundle. Then we have 
\begin{equation} \label{stability}
\int_\Sigma |A|^2 \varphi^2 \, d\mathcal{H}^n \leq \int_\Sigma \left( |\nabla \varphi|^2 + K_1 \varphi^2 \right) \, d\mathcal{H}^n
\end{equation}
for all testfunctions $\varphi \in C^\infty_c(\Sigma)$.
\end{lem}

\noindent \emph{Proof:} We test (\ref{Jacobi_flat}) with $w^{-1} \varphi^2$ and perform a partial integration. This leads to
\begin{eqnarray*}
\int_\Sigma |A|^2 \varphi^2 \, d\mathcal{H}^n & = & 2\int_\Sigma w^{-1}\varphi \nabla \varphi \nabla w \, d\mathcal{H}^n - \int_\Sigma w^{-2} |\nabla w|^2 \varphi^2 \, d\mathcal{H}^n \\
& & + \int_\Sigma w^{-1 } \Omega_{\alpha i} \nabla_i H_\alpha \varphi^2 \, d\mathcal{H}^n. 
\end{eqnarray*}
The desired estimate now follows from the Cauchy-Schwarz inequality. \qed \\

The next inequality generalizes the Simons inequality of Schoen, Simon and Yau \cite{SSY} and Ecker-Huisken \cite{Ecker:Huisken:HPC} for hypersurfaces in $\R^m$ to immersions with arbitrary co-dimension. Note that for $H=0$ we can let $\varepsilon \searrow 0$ in (\ref{Ecker_Huisken_flat}) to obtain a corresponding estimate of Smoczyk, Wang and Xin \cite{SWX} and Wang \cite{Wang:stability} for minimal immersions with flat normal bundle. 

\begin{lem} \label{Simons:inequality:flat}
Let $f:\Sigma^n \rightarrow \R^m$ be an immersion with flat normal bundle. Then we have the estimate
\begin{eqnarray}\label{Ecker_Huisken_flat}
\frac{1}{2} \Delta |A|^2 & \geq & \left( 1+ \frac{2}{n+\varepsilon} \right) |\nabla |A||^2 + h_{\alpha i j }\nabla_i \nabla_j H_\alpha \nonumber \\
& & + H_\alpha h_{\alpha ij} h_{\beta jk} h_{\beta ki} - |A|^4 - C(n,\varepsilon)|\nabla H|^2
\end{eqnarray}
for all $\varepsilon>0$.
\end{lem}

\noindent \emph{Proof:} From (\ref{Simons_flat}) we infer the estimate
\begin{equation}\label{Simons_step1}
\frac{1}{2} \Delta |A|^2  \geq  |\nabla A|^2 + h_{\alpha i j }\nabla_i \nabla_j H_\alpha \nonumber + H_\alpha h_{\alpha ij} h_{\beta jk} h_{\beta ki} - |A|^4.
\end{equation}

In any point $p_0 \in \Sigma$ where $|A|$ does not vanish, we have
$$\nabla_k |A| = |A|^{-1} \sum_{\alpha, i,j} \nabla_k h_{\alpha ij} h_{\alpha ij}.$$
Since $R^\perp_{i j \alpha \beta}=0$ we infer from the Ricci equation that we may choose our frames such that in $p_0$ all $h_{\alpha i j}$, $\alpha=n+1, \ldots, m$, are simultaneously diagonal. Hence, we obtain 
\begin{eqnarray}\label{Simons_step2}
|\nabla |A||^2  & =  & |A|^{-2}\sum_k \left( \sum_{\alpha, i} \nabla_k h_{\alpha ii} h_{\alpha ii} \right)^2 \nonumber \\
& \leq &  \sum_{\alpha, i,k} (\nabla_k h_{\alpha i i})^2 \nonumber \\
& = & \sum_{\alpha, i,k \atop i \not = k} (\nabla_k h_{\alpha i i})^2 + \sum_{\alpha, k} (\nabla_k h_{\alpha k k})^2.  
\end{eqnarray}
Moreover, we have  
\begin{eqnarray}\label{Simons_step3}
|\nabla A|^2 - |\nabla |A||^2 & \geq & \sum_{\alpha, i, j, k} (\nabla_k h_{\alpha i j})^2 - \sum_{\alpha,i,k} (\nabla_k h_{\alpha ii})^2 \nonumber \\
& = & \sum_{\alpha, i,j,k  \atop  i \not = j} (\nabla_k h_{\alpha i j})^2 \nonumber \\
 & \geq & 2 \sum_{\alpha, i,k \atop i \not = k} (\nabla_k h_{\alpha i i})^2, 
\end{eqnarray}
where the last line follows from the Codazzi equation. 

From $\nabla_k H_\alpha = \sum_i \nabla_k h_{\alpha i i}$ we infer for fixed $\alpha$ and $k$
$$(\nabla_k h_{\alpha k k})^2 = (\nabla_k H_\alpha)^2 - 2 \nabla_k H_\alpha \left( \sum_{i \atop i \not = k} \nabla_k h_{\alpha i i}\right) + \left( \sum_{i \atop i \not = k} \nabla_k h_{\alpha i i} \right)^2.$$
Applying Young's inequality and summing over $\alpha$ and $k$ leads to
\begin{equation} \label{Simons_step4}
\sum_{\alpha,k} (\nabla_k h_{\alpha k k})^2 \leq (n-1+\varepsilon) \sum_{\alpha, i,k \atop i \not = k} (\nabla_k h_{\alpha i i})^2 + \left(1+ \frac{n-1}{\varepsilon} \right) |\nabla H|^2.
\end{equation}

Combining (\ref{Simons_step1}), (\ref{Simons_step2}), (\ref{Simons_step3}) and (\ref{Simons_step4}) now gives the desired estimate (\ref{Ecker_Huisken_flat}) in all points where $|A|(p_0) \not = 0$. However, since $|A| \in W^{1, \infty}_{loc}$ with $\nabla |A|(p_0)=0$ whenever $|A|(p_0)=0$, we see that (\ref{Ecker_Huisken_flat}) must be globally true in the weak sense. \qed \\

\section{Curvature estimates}

Following Ecker-Huisken \cite{Ecker:Huisken:HPC} we define a quantity $K_2$ by
\begin{equation}\label{Definition:H2}
K_2 := \left\{\begin{array}{cl} - \left( \frac{h_{\alpha ij} \nabla_i \nabla_j H_\alpha}{|A|}\right)^- & \mbox{, if } |A|>0\\ 0 & \mbox{, if } |A|=0  \end{array}\right.,
\end{equation}
where $g^-$ denotes the negative part of the function $g$. Clearly, we have the estimate
$$K_2 \leq |\nabla \nabla H|.$$

We will now prove the following integral curvature estimate: 

\begin{thm} \label{thm:integral:curvature:estimate}
If $\Sigma^n  \subset \R^m$ is a graph with flat normal bundle, then we have
\begin{eqnarray} \label{integral:curvature:estimate}
\lefteqn{\int_\Sigma |A|^p \varphi^p \, d\mathcal{H}^n} \nonumber \\
 & \leq & C \int_\Sigma \left( |\nabla \varphi|^p + \left(|H|^p + |\nabla H|^{p/2} + K_1^{p/2} + K_2^{p/3}\right) \varphi^p  \right) \, d\mathcal{H}^n 
\end{eqnarray}
for all $p \in [4, 4+ \sqrt{8/n})$ and for all non-negative testfunctions $\varphi \in C^\infty_c(\Sigma)$, the constant $C$ depending on $n$ and $p$ only.\\ 
\end{thm} 

\noindent \emph{Proof:} We test (\ref{stability}) with $|A|^{q+1} \varphi$, where $\varphi \in C^\infty_c(\Sigma)$ is a non-negative testfunction and $q \geq 0$ is yet to be chosen, and obtain 
\begin{eqnarray}\label{int:step1}
\int_\Sigma |A|^{2q+4} \varphi^2 \, d\mathcal{H}^n & \leq & (q+1)^2 \int_\Sigma |A|^{2q}|\nabla |A||^2 \varphi^2 \, d\mathcal{H}^n \nonumber \\
& & +2(q+1) \int_\Sigma |A|^{2q+1} \varphi \nabla |A| \nabla \varphi  \, d\mathcal{H}^n \nonumber\\
& & + \int_\Sigma |A|^{2q+2} (|\nabla \varphi|^2 + K_1 \varphi^2) \, d\mathcal{H}^n.  
\end{eqnarray}

On the other hand, multiplying the Simons inequality (\ref{Ecker_Huisken_flat}) by $|A|^{2q} \varphi^2$, integrating by parts and applying Young's inequality in the form
$$|H_\alpha h_{\alpha i j}h_{\beta j k} h_{\beta k i}| \leq C(n) |H| |A|^3 \leq \varepsilon |A|^4 + \frac{C(n)}{ \varepsilon} |H|^2 |A|^2$$
leads to 
\begin{eqnarray}\label{int:step2}
\lefteqn{\left(  1+ \frac{2}{n+\varepsilon} + 2q  \right) \int_\Sigma |A|^{2q}|\nabla |A||^2 \varphi^2 \, d\mathcal{H}^n} \nonumber \\
& \leq & (1+ \varepsilon) \int_\Sigma |A|^{2q+4} \varphi^2 \, d\mathcal{H}^n  -\int_\Sigma h_{\alpha ij}\nabla_i \nabla_j H_\alpha |A|^{2q} \varphi^2 \, d\mathcal{H}^n \nonumber \\
& & +C\int_\Sigma |A|^{2q+2} |H|^2 \varphi^2 \, d\mathcal{H}^n  +C \int_\Sigma |A|^{2q} |\nabla H|^2 \varphi^2 \, d\mathcal{H}^n \nonumber \\
& & -2 \int_\Sigma |A|^{2q+1} \varphi \nabla |A| \nabla \varphi \, d\mathcal{H}^n  
\end{eqnarray}
with $C=C(n,\varepsilon)$.

Combining (\ref{int:step1}) and (\ref{int:step2}) and recalling the definition  of $K_2$ we arrive at
\begin{eqnarray} \label{int:step3}
\lefteqn{\left(1+\frac{2}{n+\varepsilon} +2q  - (1+\varepsilon)(q+1)^2 \right) \int_\Sigma |A|^{2q}|\nabla |A||^2 \varphi^2 \, d\mathcal{H}^n} \nonumber \\
& \leq & C \int_\Sigma |A|^{2q+2} (|\nabla \varphi|^2 + |H|^2 \varphi^2 + K_1 \varphi^2) \, d\mathcal{H}^n  \nonumber \\
& & + C \int_\Sigma |A|^{2q} |\nabla H|^2 \varphi^2 \, d\mathcal{H}^n +C \int_\Sigma |A|^{2q+1} K_2 \varphi^2 \, d\mathcal{H}^n \nonumber \\
& & +C \int_\Sigma |A|^{2q+1} \varphi |\nabla |A|| |\nabla \varphi| \, d\mathcal{H}^n
\end{eqnarray}
with $C=C(n,q,\varepsilon)$. We now choose $q$ such that $p=4+2q$. Then we have $q \in [0, \sqrt{2/n})$ and thus we can find $\varepsilon>0$ small enough depending on $n$ and $q$ only such that
$$1+\frac{2}{n+\varepsilon} + 2q - (1+\varepsilon)(q+1)^2 >0.$$
Hence, with this choice of $\varepsilon$ we obtain
\begin{eqnarray*}
\lefteqn{\int_\Sigma |A|^{2q} |\nabla |A||^2 \varphi^2 \, d\mathcal{H}^n}\\
& \leq  & C \int_\Sigma |A|^{2q+2} (|\nabla \varphi|^2 + |H|^2 \varphi^2 + K_1 \varphi^2) \, d\mathcal{H}^n \\
& & + C \int_\Sigma |A|^{2q} |\nabla H|^2 \varphi^2 \, d\mathcal{H}^n + C \int_\Sigma |A|^{2q+1} K_2 \varphi^2 \, d\mathcal{H}^n  \\ 
& & + C \int_\Sigma |A|^{2q+1} \varphi |\nabla |A|| |\nabla \varphi| \, d\mathcal{H}^n
\end{eqnarray*}
with $C=C(n,q)$. In view of Young's inequality and (\ref{int:step1}) this leads to
\begin{eqnarray} \label{int:step4}
\lefteqn{\int_\Sigma |A|^{2q+4} \varphi^2 \, d\mathcal{H}^n} \nonumber \\
& \leq & C \int_\Sigma |A|^{2q+2} (|\nabla \varphi|^2 + |H|^2 \varphi^2 + K_1 \varphi^2) \, d\mathcal{H}^n \nonumber \\
& & + C \int_\Sigma |A|^{2q} |\nabla H|^2 \varphi^2 \, d\mathcal{H}^n + C \int_\Sigma |A|^{2q+1} K_2 \varphi^2 \, d\mathcal{H}^n  
\end{eqnarray}
with $C=C(n,q)$.

To complete the proof we replace $\varphi$ by $\varphi^{q+2}$ in (\ref{int:step4}) and obtain
\begin{eqnarray*} 
\lefteqn{\int_\Sigma |A|^{2q+4} \varphi^{2q+4} \, d\mathcal{H}^n}\\
& \leq & C \int_\Sigma |A|^{2q+2} \varphi^{2q+2} (|\nabla \varphi|^2 + |H|^2 \varphi^2 + K_1 \varphi^2) \, d\mathcal{H}^n \\
& & + C \int_\Sigma |A|^{2q} \varphi^{2q} |\nabla H|^2 \varphi^4 \, d\mathcal{H}^n + C \int_\Sigma |A|^{2q+1} \varphi^{2q+1} K_2 \varphi^3 \, d\mathcal{H}^n 
\end{eqnarray*}
with $C=C(n,q)$. The desired inequality 
\begin{eqnarray*} 
\lefteqn{\int_\Sigma |A|^{2q+4} \varphi^{2q+4} \, d \mathcal{H}^n}\\ 
& \leq & C \int_\Sigma |\nabla \varphi|^{2q+4} \, d \mathcal{H}^n \\
& & + C \int_\Sigma \left( |H|^{2q+4} + |\nabla H|^{\frac{2q+4}{2}} + K_1^{\frac{2q+4}{2}} + K_2^{\frac{2q+4}{3}} \right) \varphi^{2q+4} \, d \mathcal{H}^n  
\end{eqnarray*}
now follows easily in view of the interpolation inequality $ab \leq \gamma a^s + \gamma^{-\frac{t}{s}} b^t$ for all $a,b \geq 0$, $\gamma>0$ and $s,t>1$ with $\frac{1}{s} + \frac{1}{t}=1$. \qed\\

Denote by $B_R = B_R(p) \subset \R^m$ the closed ball of radius $R>0$ with center $p \in \Sigma$. In order to obtain a $sup$ curvature estimate we need the following mean value inequality. The proof is similar to \cite[Theorem 8.17]{Gilbarg:Trudinger}, however we assume less regularity on the coefficients of (\ref{subsolution}). For the convenience of the reader we sketch the argument below.  

\begin{lem} \label{meanvalue}
Let $\Sigma^n \subset \R^m$ be an arbitrary graph, and suppose that $u$ is a non-negative solution of
\begin{equation} \label{subsolution}
\Delta u + Qu \geq g \quad \mbox{on } \Sigma, 
\end{equation}
where $Q \in L^{q/2}(\Sigma)$ and $g \in L^{p/2}(\Sigma)$ with $q,p > n$. If $\Sigma \cap B_{2R}\subset \subset \Sigma$ then we have the estimate
\begin{equation} \label{supestimate}
\sup_{\Sigma \cap B_R} u \leq C\left( R^{-n/2} \|u\|_{L^2(\Sigma \cap B_{2R})} + k(R) \right),
\end{equation}
where
\begin{equation} \label{definitionk}
k(R) = R^{2(1-n/p)}\|g\|_{L^{p/2}(\Sigma \cap B_{2R})}, 
\end{equation}
the constant $C$ depending on $n$, $q$, $p$, $R^{2(1-n/q)}\|Q\|_{L^{q/2}(\Sigma \cap B_{2R})}$, $R\sup_{\Sigma \cap B_{2R}} |H|$ and $R^{-n} \mathcal{H}^n(\Sigma \cap B_{2R})$.
\end{lem}

\noindent \emph{Proof:} First, note that by scaling $\R^m \rightarrow \R^m, p \mapsto R p$ it suffices to consider the case $R=1$. 

We now put $v=u+k$, where $k=\|g\|_{L^{p/2}(\Sigma \cap B_2)}$, and let $\eta \in C^\infty_c(\Sigma)$ be a non-negative function supported in $\Sigma \cap B_2$. For $\beta \geq 1$ we multiply (\ref{subsolution}) with $v^\beta \eta^2$ and perform a partial integration. This leads to
\begin{eqnarray} \label{lemma_step1}
\beta \int_\Sigma v^{\beta-1} |\nabla u|^2 \eta^2 \, d\mathcal{H}^n & \leq & -2 \int_\Sigma v^\beta \eta \nabla \eta \nabla u \, d\mathcal{H}^n \nonumber \\
& & + \int_\Sigma (Qu-g) v^\beta \eta^2 \, d\mathcal{H}^n.
\end{eqnarray}
Using Young's inequality we find
\begin{equation} \label{lemma_step2}
|2v^\beta \eta \nabla \eta \nabla u| \leq \frac{\beta}{2} v^{\beta-1} |\nabla u|^2 \eta^2 + \frac{2}{\beta} v^{\beta+1}|\nabla \eta|^2.
\end{equation}
Furthermore, since $v \geq \max(u,k)$ we have
\begin{equation} \label{lemma_step3}
| (Qu -g) v^\beta \eta^2 | \leq v^{\beta +1} \eta^2 \left(|Q| + \frac{|g|}{k} \right),
\end{equation}
where $\frac{|g|}{k}$ is to be considered $0$ in case $k=0$.
Combining (\ref{lemma_step1}), (\ref{lemma_step2}) and (\ref{lemma_step3}) yields
\begin{eqnarray*}
\int_\Sigma v^{\beta-1}  |\nabla u|^2 \eta^2 \, d\mathcal{H}^n & \leq & \frac{4}{\beta^2} \int_\Sigma v^{\beta+1} |\nabla \eta|^2  \, d\mathcal{H}^n \\
& & + \frac{2}{\beta} \int_\Sigma  v^{\beta+1} \eta^2 \left( |Q| + \frac{|g|}{k} \right) \, d\mathcal{H}^n.
\end{eqnarray*}
Hence, abbreveating $w=v^{\frac{\beta+1}{2}}$ we arrive at the estimate
\begin{eqnarray} \label{lemma_step4}
\int_\Sigma  |\nabla w|^2 \eta^2 \, d\mathcal{H}^n & \leq & 4 \int_\Sigma  w^2 |\nabla \eta|^2  \, d\mathcal{H}^n \nonumber \\
& & + 2 \beta \int_\Sigma  w^2 \eta^2 \left( |Q| + \frac{|g|}{k} \right) \, d\mathcal{H}^n.
\end{eqnarray}

Next, we apply the Sobolev-inequality of Michael-Simon \cite{Michael:Simon} followed by Hölder's inequality to obtain
\begin{eqnarray*}
\lefteqn{\left( \int_\Sigma (\eta w)^{2\chi} \, d\mathcal{H}^n \right)^{\frac{1}{\chi}}} \nonumber \\
& \leq & C \int_\Sigma (|\nabla \eta|^2 w^2 +  \eta^2 |\nabla w|^2 + \eta^2 w^2 |H|^2) \, d\mathcal{H}^n,
\end{eqnarray*} 
where $\chi=\frac{\hat{n}}{\hat{n}-2}$ with $\hat{n}=n $ for $n \geq 3$ and $2 < \hat{n} < \min\{q,p\}$ for $n=2$, respectively, and where $C$ is a constant depending on $\hat{n}$ and $\mathcal{H}^n(\Sigma \cap B_2)$. Combining this with (\ref{lemma_step4}) leads to
\begin{eqnarray} \label{lemma_step5}
\left( \int_\Sigma (\eta w)^{2\chi}  \, d\mathcal{H}^n \right)^{\frac{1}{\chi}} & \leq & C \int_\Sigma w^2 (\eta^2 + |\nabla \eta|^2 )\, d\mathcal{H}^n \nonumber \\
& & + C \beta \int_\Sigma  w^2 \eta^2 \left( |Q| + \frac{|g|}{k} \right) \, d\mathcal{H}^n,
\end{eqnarray}
the constant $C$ now depending additionally on $\sup_{\Sigma \cap B_2} |H|$.

Next we use interpolation inequalities for $L^p$-spaces, cf. \cite[Section 7.1]{Gilbarg:Trudinger}, and obtain
\begin{eqnarray}\label{lemma_step6}
\lefteqn{\int_\Sigma w^2 \eta^2 |Q| \, d\mathcal{H}^n} \nonumber \\
& \leq & \left( \int_\Sigma (w\eta)^{\frac{2q}{q-2}}\, d\mathcal{H}^n \right)^{\frac{q-2}{q}} \left( \int_{\Sigma \cap B_2} |Q|^{\frac{q}{2}} \, d\mathcal{H}^n \right)^{\frac{2}{q}} \\
& \leq & \left[ \varepsilon \left( \int_\Sigma (w\eta)^{2\chi} \, d\mathcal{H}^n \right)^{\frac{1}{2\chi}} + \varepsilon^{-\mu} \left( \int_\Sigma w^2 \eta^2 \, d\mathcal{H}^n \right)^{\frac{1}{2}} \right]^2 \|Q\|_{L^{q/2}(\Sigma \cap B_2)} \nonumber
\end{eqnarray}
for all $\varepsilon >0$ with $\mu=\frac{\hat{n}}{q-\hat{n}}>0$. Similarly, we have
\begin{eqnarray}\label{lemma_step7}
\lefteqn{\int_\Sigma w^2 \eta^2 \frac{|g|}{k} \, d\mathcal{H}^n} \nonumber\\
& \leq &  \left[ \varepsilon \left( \int_\Sigma (w\eta)^{2 \chi} \, d\mathcal{H}^n \right)^{\frac{1}{2\chi}} + \varepsilon^{-\tilde{\mu}} \left( \int_\Sigma w^2 \eta^2 \, d\mathcal{H}^n \right)^{\frac{1}{2}}\right]^2 
\end{eqnarray}
with $\tilde{\mu} = \frac{\hat{n}}{p-\hat{n}}>0$. Hence, using (\ref{lemma_step6}), (\ref{lemma_step7}) with $\varepsilon \sim [\beta (\|Q\|_{L^{q/2}(\Sigma \cap B_2)} +1)]^{-1/2}$ in (\ref{lemma_step5}) we finally arrive at 
\begin{equation} \label{lemma_step8}
\left( \int_\Sigma (\eta w)^{2\chi}\, d\mathcal{H}^n \right)^{\frac{1}{\chi}} \leq C \beta^\alpha \int_\Sigma w^2(\eta^2 + |\nabla \eta|^2 )\, d\mathcal{H}^n
\end{equation}
with $C$ depending on $n$, $q$, $p$, $\sup_{\Sigma \cap B_2} |H|$, $\mathcal{H}^n(\Sigma \cap B_2)$ and $\|Q\|_{L^{q/2}(\Sigma \cap B_2)}$, and $\alpha=\alpha(n,p,q)>1$. 

From here we can employ Moser's iteration technique \cite{Moser} in a manner similar to \cite{Ecker:Huisken:HPC} and \cite[Section 4]{winklmann:manuscripta}. Put $\gamma := \beta +1 \geq 2$ such that $w^2=v^\gamma$. Let $\rho, \rho'$ be radii satisfying $1 \leq \rho' \leq \rho \leq 2$ and let $\eta \in C^\infty_c(\Sigma)$ to be a cut-off function with $0 \leq \eta \leq 1$, $\eta=1$ in $\Sigma \cap B_{\rho'}$, $\mbox{supp} (\eta) \subset \Sigma \cap B_\rho$, and $|\nabla \eta| \leq \frac{C}{\rho-\rho'}$. Then we infer from (\ref{lemma_step8}) the estimate
\begin{eqnarray}\label{lemma_step9}
\left( \int_{\Sigma \cap B_{\rho'}} v^{\chi \gamma} \, d\mathcal{H}^n  \right)^{\frac{1}{\chi \gamma}} \leq \frac{C^{\frac{1}{\gamma}} \gamma^{\frac{\alpha}{\gamma}}}{(\rho-\rho')^{\frac{2}{\gamma}}} \left( \int_{\Sigma \cap B_{\rho}} v^{\gamma} \, d\mathcal{H}^n \right)^{\frac{1}{\gamma}} 
\end{eqnarray} 
with a constant $C$ depending on $n$, $q$, $p$, $\sup_{\Sigma \cap B_2} |H|$, $\mathcal{H}^n(\Sigma \cap B_2)$ and $\|Q\|_{L^{q/2}(\Sigma \cap B_2)}$ only.  Now, let 
$$\rho_k = 1+2^{-k}, \quad \rho_k' = \rho_{k+1}, \quad \gamma_k = 2 \chi^k\quad \mbox{for }k=0,1,2,\ldots.$$ 
Replacing $\rho$, $\rho'$ and $\gamma$ in (\ref{lemma_step9}) by $\rho_k$, $\rho_k'$ and $\gamma_k$ and iterating the resulting inequalities as $k \rightarrow \infty$, we obtain the estimate
$$\sup_{\Sigma \cap B_1} v \leq C\left( \int_{\Sigma \cap B_2} v^2 \, d\mathcal{H}^n \right)^{\frac{1}{2}}$$
with $C$ depending on the same data as before. Recalling that $v=u+k$, this gives the desired result. \qed \\

Now we are ready to prove our main result. 

\begin{thm} \label{main:thm}
Let $\Sigma^n \subset \R^m$, $2 \leq n \leq 5$,  be a graph with flat normal bundle, and suppose that $\Sigma \cap B_{4R} \subset \subset \Sigma$ with 
$$\mathcal{H}^n(\Sigma \cap B_{4R}) \leq K R^n. $$ 
Then we have the estimate
\begin{equation}\label{sup_estimate}
\sup_{\Sigma \cap B_R} |A|^2 \leq \frac{C}{R^2} 
\end{equation} 
with a constant $C$ depending on $n$, $K$, $R \sup_{\Sigma \cap B_{4R}} |H|$, $R^2 \sup_{\Sigma \cap B_{4R}} (|\nabla H| +  K_1)$ and $R^3 \sup_{\Sigma \cap B_{4R}} K_2$.\\
\end{thm}

\noindent \emph{Proof:} In view of the Simons identity (\ref{Simons_flat}) and the estimate
$$|H_\alpha h_{\alpha i j}h_{\beta j k} h_{\beta k i}| \leq C(n)|A|^4$$ 
we infer 
$$\Delta |A|^2 + C(n)|A|^4 \geq -2K_2|A|.$$

Furthermore, since $2 \leq n \leq 5$, we can apply Theorem \ref{thm:integral:curvature:estimate} with a suitable cut-off function as before to obtain     
$$\int_{\Sigma \cap B_{2R}} |A|^q \, d\mathcal{H}^n \leq CR^{n-q}$$
for some $q > \max\{n,4\}$, with a constant $C$ depending on $n$, $K$, $R \sup_{\Sigma \cap B_{4R}} |H|$, $R^2 \sup_{\Sigma \cap B_{4R}} (|\nabla H| +  K_1)$ and $R^3 \sup_{\Sigma \cap B_{4R}} K_2$.

Hence, applying Lemma \ref{meanvalue} with $u=|A|^2$, $Q=C(n)|A|^2$, $g=-2K_2|A|$ and $p=2q$, the desired estimate follows easily. \qed\\ 

In case $H=0$ the constant in Theorem \ref{main:thm} is independent of $R$. Therefore, letting $R \rightarrow \infty$ in (\ref{sup_estimate}) we obtain the Bernstein result of Smoczyk, Wang and Xin \cite{SWX} and Wang \cite{Wang:stability}:

\begin{cor}
Suppose that $\Sigma = \{(x,\psi(x)) : x \in \R^n\} \subset \R^m$, $2 \leq n \leq 5$, is an entire minimal graph with flat normal bundle. If 
$$\mathcal{H}^n(\Sigma \cap B_R(p) ) \leq KR^n$$
for some point $p \in \Sigma$ and some sequence $R \rightarrow \infty$ with a constant $K$ independent of $R$, 
then $\psi$ is an affine linear function.  
\end{cor}

\vspace{4ex}

\noindent Steffen Fröhlich
 
\noindent Technische Universität Darmstadt

\noindent FB Mathematik

\noindent Schloßgartenstraße 7

\noindent 64289 Darmstadt, Germany

\noindent sfroehlich@mathematik.tu-darmstadt.de\\

\noindent Sven Winklmann

\noindent Centro di Ricerca Matematica Ennio De Giorgi 

\noindent Scuola Normale Superiore di Pisa

\noindent Piazza dei Cavalieri 3

\noindent 56100 Pisa, Italy

\noindent s.winklmann@sns.it

\end{document}